\documentclass[a4paper]{article}
\usepackage{amsmath,amssymb,amsfonts}
\usepackage{graphicx} 
\usepackage{float} 
\usepackage{color}
\usepackage{lscape}

\definecolor{Azul}{cmyk}{1,0.50,0,0.3}
\definecolor{Magenta}{cmyk}{0,1,0,0}

\newtheorem{theorem}{Theorem}[section]

\newtheorem{lemma}{Lemma}[section]

\title{The non-trivial zeros of Riemann's zeta-function}
\author{Jailton C. Ferreira}

\date{ }
\begin{document}
\maketitle \pagenumbering{arabic}

\begin{abstract}
A proof of the Riemann hypothesis using the reflection principle
$\zeta( \overline s) = \overline {\zeta(s)}$ is presented.
\end{abstract}

\section{Introduction} \label{sec-1}

\hspace{22pt} The Riemann zeta-function $ \zeta (s) $ can be
defined by either of two following formulae \cite{Titchmarsh}

\begin{equation}\label{um-1}
\zeta(s) = \sum_{n=1}^{\infty} \frac{1}{n^s}
\end{equation}

where $n \in \mathbb{N}$, and
\begin{equation}\label{um-2}
\zeta(s) = \prod_p \left( 1 - \frac{1}{p^s} \right)^{-1}
\end{equation}

where $p$ runs through all primes and $ \ \Re(s) > 1$. By analytic
continuation $\zeta(s)$ is defined over whole $ \mathbb{C}$.

\hspace{22pt} The relationship between $\zeta(s)$ and $\zeta(1 -
s)$
\begin{equation}\label{um-3}
\zeta(s) = 2^s \pi^{s-1} \sin \left( \frac{\pi s}{2} \right)
\Gamma( 1 - s) \zeta(1 - s)
\end{equation}

is known as the functional equation of the zeta-function. From the
functional equation it follows that $\zeta(s)$ has zeros at $s =
-2, -4, -6, \ldots $ . These zeros are traditionally called
trivial zeros of $\zeta(s)$; the zeros of $\zeta(s)$ with $\Im(z)
\ne 0$ are called non-trivial zeros. From the equation
\eqref{um-2}, which is known as Euler's product, it was deduced
that $\zeta(s)$ has no zeros for $ \Re(s) > 1 $. The functional
equation implies that there are no non-trivial zeros with $ \Re(s)
< 0 $. It was deduced that there are no zeros for $\Re(s) = 0$ and
$\Re(s) = 1$. Therefore all non-trivial zeros are in the
\textit{critical strip} specified by $ 0 < \Re(s) < 1$.

\hspace{22pt} In 1859 Riemann published the paper \textit{Ueber
die Anzahl der Primzahlen unter einer gegebenen Gr\"osse}. A
translation of the paper is found in \cite{Edwards}. In the paper
Riemann considers ``very likely" that all the non-trivial zeros of
$\zeta(s)$ have real part equal to $ \frac{1}{2} $. The statement
\begin{center}
\textit{The non-trivial zeros of $\zeta(s)$ have real part equal
to $ \frac{1}{2} $.}
\end{center}
is known as the Riemann hypothesis.

\section{Theorem} \label{sec-2}

\begin{lemma} \label{lemma-1}
If $z$ is such that $ 0 < \Re(z) < 1$,
$\Im(z) \ne 0$ and

\begin{equation}\label{lema1-1}
\zeta(z) =0
\end{equation}
then
\begin{equation}\label{lema1-2}
 \zeta(1 - z)=0
\end{equation}
\end{lemma}
\textit{Proof:}

\hspace{22pt} Let us examine the product
\begin{equation}\label{lema1-2}
 2^s \pi^{s-1} \sin \left( \frac{\pi s}{2} \right) \Gamma( 1 - s) 
\end{equation}
for  $ 0 \le \Re(s) \le 1$.

\hspace{22pt} In the Figures 1 to 10 the curves in blue, red and brown are associated with $\Re(s)=0$, $\Re(s)=0.5$ and $\Re(s)=1$, respectively.

\hspace{22pt} The Figures 1 and 2 show $ 2^s$. The function  $ 2^s$ is not equal to 0 in the range $ 0 \le \Re(s) \le 1$, The figures illustrated this.

\begin{figure}
\begin{center}
\includegraphics{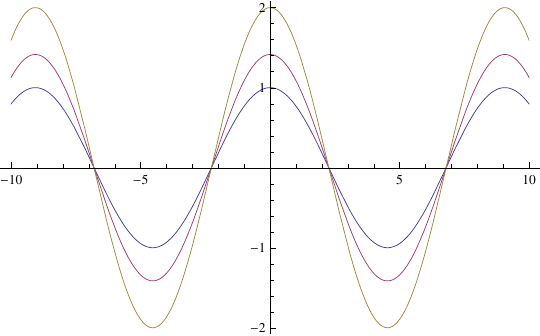}
\end{center}
\caption{Real part of $2^s$ versus imaginary part of $s$}
\label{figure1}
\end{figure} 

\begin{figure}
\begin{center}
\includegraphics{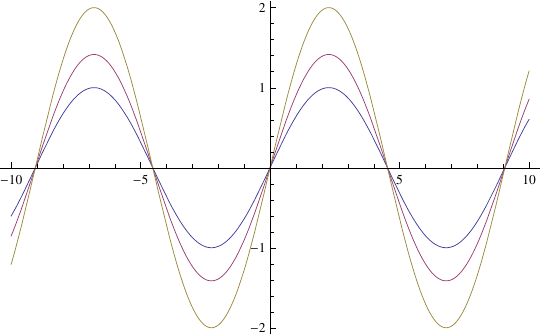}
\end{center}
\caption{Imaginary part of $2^s$ versus imaginary part of $s$}
\label{figure2}
\end{figure} 

\hspace{22pt} The Figures 3 and 4 show $ \pi^{s-1}$. The function  $ \pi^{s-1}$ is not equal to 0 in the range $ 0 \le \Re(s) \le 1$, The figures illustrated this.

\begin{figure}
\begin{center}
\includegraphics{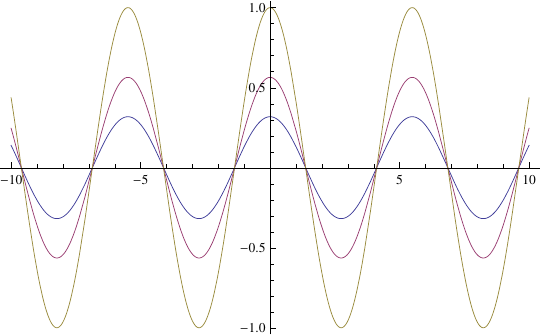}
\end{center}
\caption{Real part of $\pi^{s-1}$ versus imaginary part of $s$}
\label{figure3}
\end{figure} 

\begin{figure}
\begin{center}
\includegraphics{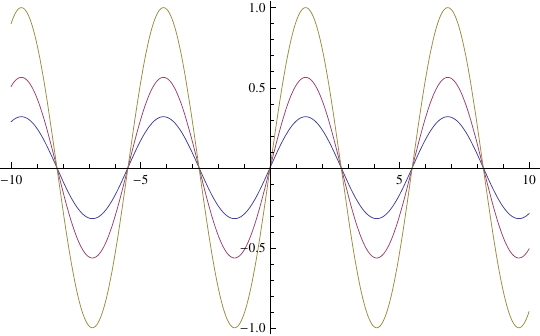}
\end{center}
\caption{Imaginary part of $\pi^{s-1}$ versus imaginary part of $s$}
\label{figure4}
\end{figure} 

\hspace{22pt} The Figures 5 to 8 show $ \sin \left( \frac{\pi s}{2} \right)$. The function  $ \sin \left( \frac{\pi s}{2} \right)$ is not equal to 0 in the range $ 0 \le \Re(s) \le 1$, except for $\Im(s) = 0$. The figures illustrated this.

\begin{figure}
\begin{center}
\includegraphics{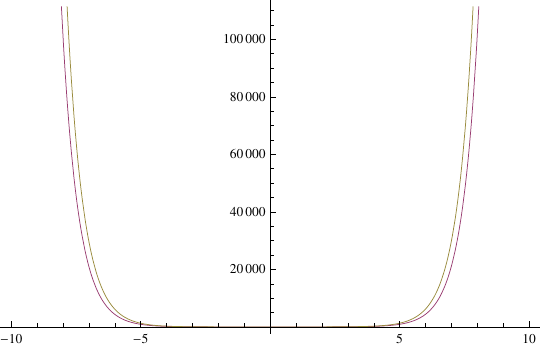}
\end{center}
\caption{Real part of $ \sin \left( \frac{\pi s}{2} \right)$ versus imaginary part of $s$}
\label{figure5}
\end{figure} 

\begin{figure}
\begin{center}
\includegraphics{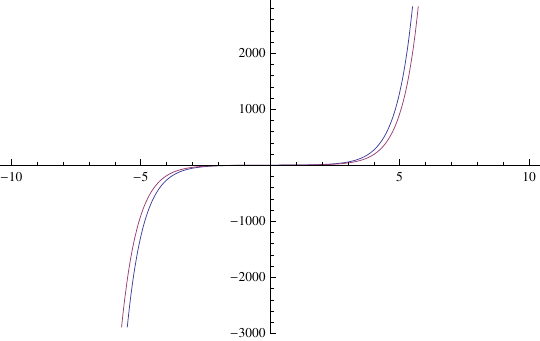}
\end{center}
\caption{Imaginary part of $ \sin \left( \frac{\pi s}{2} \right)$ versus imaginary part of $s$}
\label{figure6}
\end{figure} 

\begin{figure}
\begin{center}
\includegraphics{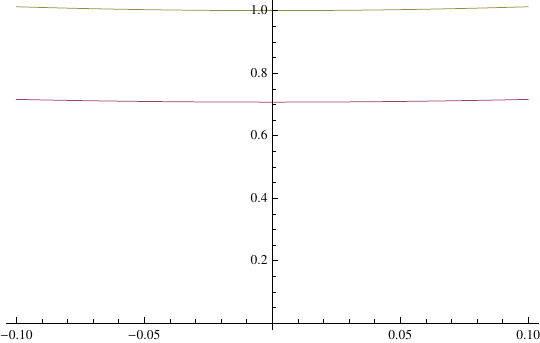}
\end{center}
\caption{Real part of $ \sin \left( \frac{\pi s}{2} \right)$ versus imaginary part of $s$}
\label{figure7}
\end{figure} 

\begin{figure}
\begin{center}
\includegraphics{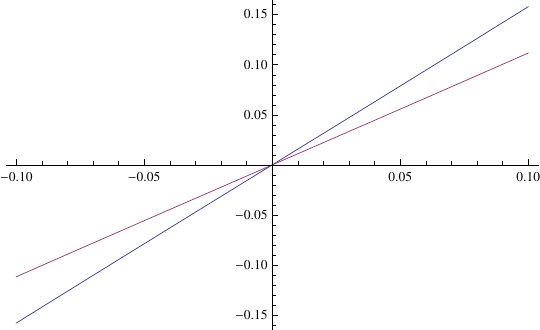}
\end{center}
\caption{Imaginary part of $ \sin \left( \frac{\pi s}{2} \right)$ versus imaginary part of $s$}
\label{figure8}
\end{figure} 

\hspace{22pt} The Figures 9 and 10 show $ \Gamma(1 - s)$. It is known that there is no complex number $s$ for which $ \Gamma(s)=0 $.

\begin{figure}
\begin{center}
\includegraphics{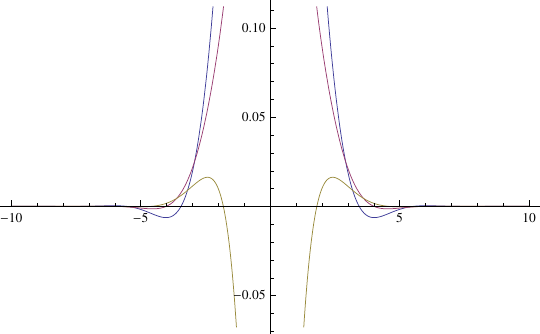}
\end{center}
\caption{Real part of $ \Gamma(1 - s)$ versus imaginary part of $s$}
\label{figure9}
\end{figure} 

\begin{figure}
\begin{center}
\includegraphics{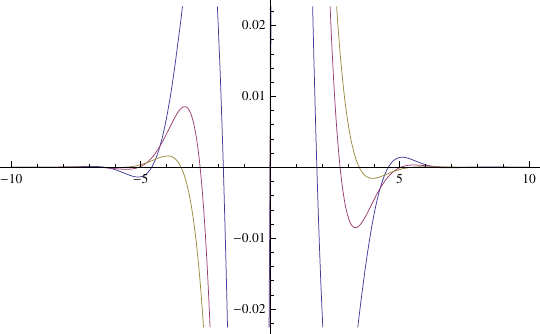}
\end{center}
\caption{Imaginary part of $ \Gamma(1 - s)$ versus imaginary part of $s$}
\label{figure10}
\end{figure} 

\hspace{22pt} The product given in \eqref{lema1-2} is different of zero in the range $ 0 \le \Re(s) \le 1$. Considering this we can conclude from \eqref{um-3} that if $s$ is such that $ 0 < \Re(s) < 1$, $\Im(s) \ne 0$ and $\zeta(s) =0$ then $ \zeta(1 - s)=0$.

\hspace{400pt} $\Box$

\begin{lemma} \label{lemma-2}
If
\begin{equation}\label{lema-1}
\sum_{n=1}^{\infty} \frac{1}{n^{2
\delta}}\frac{(-1)^{n-1}}{n^{1-z}} = \sum_{n=1}^{\infty}
\frac{(-1)^{n-1}}{n^{1-z}}= 0
\end{equation}
where
\begin{equation}\label{lema-2}
\Re(z) = \frac{1}{2} + \delta
\end{equation}
\begin{equation}\label{lema-3d}
\Im(z) = t \ne 0 \
\end{equation}
and
\begin{equation}\label{lema-4}
-\frac{1}{2} < \delta < \frac{1}{2}
\end{equation}
then
\begin{equation}\label{lema-5}
\delta = 0
\end{equation}
\end{lemma}

\textit{Proof:}

\hspace{22pt} Let be
\begin{equation}\label{lema-6}
\frac{(-1)^{n-1}}{n^{1-z}} = b_n + i c_n \qquad \textrm{for} \quad
n=1,2,3, \ldots
\end{equation}

From \eqref{lema-1} and \eqref{lema-6} we obtain
\begin{equation}\label{lema-7}
\sum_{n=1}^{\infty} \frac{1}{n^{2\delta}}b_n = \sum_{n=1}^{\infty}
b_n = 0
\end{equation}
and
\begin{equation}\label{lema-8}
\sum_{n=1}^{\infty} \frac{1}{n^{2\delta}}c_n = \sum_{n=1}^{\infty}
c_n = 0
\end{equation}

Since that
\begin{equation}\label{lema-9}
\sum_{n=1}^{\infty} \frac{(-1)^{n-1}}{n^{1-z}} =
\sum_{n=1}^{\infty} \Big( \frac{1}{(2n-1)^{1-z}}- \frac{1}{(2n)^{1-z}} \Big)
\end{equation}
and
\begin{equation}\label{lema-10}
\frac{1}{(2n)^{1-z}} = \frac{(2n)^{it}}{(2n)^{\frac{1}{2} -
\delta}} = \frac{1}{(2n)^{\frac{1}{2} - \delta}} \left[ \cos(t \log(2n)) +
i \sin(t \log(2n)) \right]
\end{equation}
and
\begin{equation}\label{lema-11}
\frac{1}{(2n-1)^{1-z}} = \frac{(2n-1)^{it}}{(2n-1)^{\frac{1}{2} -
\delta}} = \frac{1}{(2n-1)^{\frac{1}{2} - \delta}} \left[ \cos(t
\log(2n-1)) + i \sin(t \log(2n-1)) \right]
\end{equation}
we obtain
\begin{equation}\label{lema-12}
\sum_{n=1}^{\infty}b_n = \sum_{n=1}^{\infty} \Big( \frac{\cos(t
\log(2n-1))} {(2n-1)^{\frac{1}{2} - \delta}} - \frac{\cos(t
\log(2n))} {(2n)^{\frac{1}{2} - \delta}} \Big)
\end{equation}
and
\begin{equation}\label{lema-13}
\sum_{n=1}^{\infty}c_n = \sum_{n=1}^{\infty}  \Big( \frac{\sin(t
\log(2n-1))} {(2n-1)^{\frac{1}{2} - \delta}} - \frac{\sin(t
\log(2n))} {(2n)^{\frac{1}{2} - \delta}} \Big)
\end{equation}

Defining the functions
\begin{equation}\label{lema-14}
u(n) = \frac{\cos(t \log(2n-1))} {(2n-1)^{\frac{1}{2}}}
\end{equation}
and
\begin{equation}\label{lema-15}
v(n) = \frac{\cos(t \log(2n))} {(2n)^{\frac{1}{2}}}
\end{equation}

substituting \eqref{lema-14} and \eqref{lema-15} into
\eqref{lema-12} and using \eqref{lema-7} we have
\begin{equation}\label{lema-16}
\sum_{n=1}^{\infty} \Big( (2n-1)^{\delta} u(n) - (2n)^{\delta}
v(n) \Big) = 0
\end{equation}
and
\begin{equation}\label{lema-17}
\sum_{n=1}^{\infty} n^{-2 \delta} \Big( (2n-1)^{\delta} u(n) -
(2n)^{\delta} v(n) \Big) = 0
\end{equation}

\hspace{22pt} The $k$-th part of \eqref{lema-16} is
\begin{eqnarray}\label{lema-18}
&   & (2k-1)^{\delta} u(k) - (2k)^{\delta} v(k) = \nonumber \\
& - & \sum_{n=1}^{k-1} \Big( (2n-1)^{\delta} u(n) - (2n)^{\delta}
v(n) \Big) - \sum_{n=k+1}^{\infty} \Big( (2n-1)^{\delta} u(n) -
(2n)^{\delta} v(n) \Big)
\end{eqnarray}

where $k > 1$. From \eqref{lema-18} we can obtain the $k$-th part
of \eqref{lema-17}
\begin{eqnarray}\label{lema-19}
&   & k^{-2 \delta} \Big( (2k-1)^{\delta} u(k) - (2k)^{\delta} v(k) \Big) = \nonumber \\
& - & \sum_{n=1}^{k-1} k^{-2 \delta} \Big( (2n-1)^{\delta} u(n) -
(2n)^{\delta} v(n) \Big) - \sum_{n=k+1}^{\infty} k^{-2 \delta}
\Big( (2n-1)^{\delta} u(n) - (2n)^{\delta} v(n) \Big)
\end{eqnarray}

The left side of \eqref{lema-19} obtained from \eqref{lema-17} is
\begin{eqnarray}\label{lema-20}
&   & k^{-2 \delta} \Big( (2k-1)^{\delta} u(k) - (2k)^{\delta} v(k) \Big) = \nonumber \\
& - & \sum_{n=1}^{k-1} n^{-2 \delta} \Big( (2n-1)^{\delta} u(n) -
(2n)^{\delta} v(n) \Big) - \sum_{n=k+1}^{\infty} n^{-2 \delta}
\Big( (2n-1)^{\delta} u(n) - (2n)^{\delta} v(n) \Big)
\end{eqnarray}

Comparing \eqref{lema-19} with \eqref{lema-20} we conclude
\begin{eqnarray}\label{lema-21}
&   & \sum_{n=1}^{k-1} k^{-2 \delta} \Big( (2n-1)^{\delta} u(n) -
(2n)^{\delta} v(n) \Big) + \sum_{n=k+1}^{\infty} k^{-2 \delta}
\Big( (2n-1)^{\delta} u(n) - (2n)^{\delta} v(n) \Big) = \nonumber \\
&   & \sum_{n=1}^{k-1} n^{-2 \delta} \Big( (2n-1)^{\delta} u(n) -
(2n)^{\delta} v(n) \Big) + \sum_{n=k+1}^{\infty} n^{-2 \delta}
\Big( (2n-1)^{\delta} u(n) - (2n)^{\delta} v(n) \Big)
\end{eqnarray}
for all $k > 1$.

\hspace{22pt} Let be $\delta > 0$. Rearranging \eqref{lema-21} we have
\begin{eqnarray}\label{lema-22}
&   & \sum_{n=1}^{k-1} \Big( (2n-1)^{\delta} u(n) -
(2n)^{\delta} v(n) \Big) =- \sum_{n=k+1}^{\infty} 
\Big( (2n-1)^{\delta} u(n) - (2n)^{\delta} v(n) \Big) + \nonumber \\
&   & \sum_{n=1}^{k-1} \big( \frac{k}{n} \big) ^{2 \delta} \Big( (2n-1)^{\delta} u(n) -
(2n)^{\delta} v(n) \Big) + \sum_{n=k+1}^{\infty} \big( \frac{k}{n} \big) ^{2 \delta}
\Big( (2n-1)^{\delta} u(n) - (2n)^{\delta} v(n) \Big)
\end{eqnarray}

Considering \eqref{lema-16} the limit of the right side of \eqref{lema-22} when $k$ tends to infinite is

\begin{multline}\label{lema-23}
\lim_{k \rightarrow \infty}\Bigg[- \sum_{n=k+1}^{\infty} 
\Big( (2n-1)^{\delta} u(n) - (2n)^{\delta} v(n) \Big)  + \sum_{n=1}^{k-1} \big( \frac{k}{n} \big) ^{2 \delta} \Big( (2n-1)^{\delta} u(n) -
(2n)^{\delta} v(n) \Big)
 \\
+ \sum_{n=k+1}^{\infty} \big( \frac{k}{n} \big) ^{2 \delta}
\Big( (2n-1)^{\delta} u(n) - (2n)^{\delta} v(n) \Big)\Bigg]=0
\end{multline}

To take the limit of a sum of functions all we need to do is take the limit of the individual functions and sum them. We have that
\begin{equation}\label{lema-24}
\lim_{k \rightarrow \infty}\Bigg[ \sum_{n=k+1}^{\infty} \Big( (2n-1)^{\delta} u(n) - (2n)^{\delta} v(n) \Big) \Bigg]= \lim_{n \rightarrow \infty}\Big( (2n-1)^{\delta} u(n) - (2n)^{\delta} v(n) \Big)
\end{equation}
and
\begin{equation}\label{lema-25}
 \lim_{n \rightarrow \infty}  \frac{ (2n-1)^{\delta}} {(2n-1)^{\frac{1}{2}}}= \lim_{n \rightarrow \infty}  \frac{ (2n)^{\delta}} {(2n)^{\frac{1}{2}}}=0
\end{equation}
Considering \eqref{lema-14}, \eqref{lema-15}, \eqref{lema-24} and \eqref{lema-25} we obtain
\begin{equation}\label{lema-26}
\lim_{k \rightarrow \infty}\Bigg[ \sum_{n=k+1}^{\infty} \Big( (2n-1)^{\delta} u(n) - (2n)^{\delta} v(n) \Big) \Bigg]= 0
\end{equation}
and the equation \eqref{lema-23} becomes
\begin{multline}\label{lema-27}
\lim_{k \rightarrow \infty}\Bigg[ \sum_{n=1}^{k-1} \big( \frac{k}{n} \big) ^{2 \delta} \Big( (2n-1)^{\delta} u(n) -
(2n)^{\delta} v(n) \Big)
 \\
+ \sum_{n=k+1}^{\infty} \big( \frac{k}{n} \big) ^{2 \delta}
\Big( (2n-1)^{\delta} u(n) - (2n)^{\delta} v(n) \Big)\Bigg]=0
\end{multline}
We have that
\begin{equation}\label{lema-28}
\lim_{k \rightarrow \infty}\Bigg[ \sum_{n=k+1}^{\infty}  \big( \frac{k}{n} \big) ^{2 \delta}\Big( (2n-1)^{\delta} u(n) - (2n)^{\delta} v(n) \Big) \Bigg]= \lim_{n \rightarrow \infty} \big( \frac{n-1}{n} \big) ^{2 \delta}\Big( (2n-1)^{\delta} u(n) - (2n)^{\delta} v(n) \Big)
\end{equation}
therefore
\begin{equation}\label{lema-29}
\lim_{k \rightarrow \infty}\Bigg[ \sum_{n=k+1}^{\infty}  \big( \frac{k}{n} \big) ^{2 \delta}\Big( (2n-1)^{\delta} u(n) - (2n)^{\delta} v(n) \Big) \Bigg]= 0
\end{equation}
and \eqref{lema-27} becomes
\begin{equation}\label{lema-30}
\lim_{k \rightarrow \infty}\Bigg[ \sum_{n=1}^{k-1} \big( \frac{k}{n} \big) ^{2 \delta} \Big( (2n-1)^{\delta} u(n) -
(2n)^{\delta} v(n) \Big) \Bigg]=0
\end{equation}
or
\begin{equation}\label{lema-31}
\lim_{k \rightarrow \infty}\Bigg[ \sum_{n=1}^{k-1} \big( \frac{k}{n} \big) ^{2 \delta} \Big( (2n-1)^{\delta}\frac{\cos(t \log(2n-1))} {(2n-1)^{\frac{1}{2}}} -(2n)^{\delta}  \frac{\cos(t \log(2n))} {(2n)^{\frac{1}{2}}} \Big) \Bigg]=0
\end{equation}
The limit of the first addend of the sum in \eqref{lema-31} is
\begin{equation}\label{lema-32}
\lim_{k \rightarrow \infty}\Bigg[k^{2 \delta} \Big(\cos(t \times 0) -2^{\delta}  \frac{\cos(t \log(2))} {2^{\frac{1}{2}}} \Big) \Bigg]
\end{equation}
and
\begin{equation}\label{lema-33}
1 - 2^{\delta}  \frac{\cos(t \log(2))} {2^{\frac{1}{2}}}>0
\end{equation}
so
\begin{equation}\label{lema-34}
\lim_{k \rightarrow \infty}\Bigg[k^{2 \delta} \Big(1 -2^{\delta}  \frac{\cos(t \log(2))} {2^{\frac{1}{2}}} \Big) \Bigg]=+\infty
\end{equation}
The equation \eqref{lema-34} implies that \eqref{lema-31} does not hold.

\hspace{22pt} Considering the symmetry conditions of the zeros of Riemann's zeta-function, we conclude that $\delta$ can not belongs to $(- \frac{1}{2}, 0 )$. Hence $\delta = 0$.

\hspace{400pt} $\Box$

\begin{theorem} \label{theorem}
The non-trivial zeros of $\zeta(s)$ have real part equal to $
\frac{1}{2} $.
\end{theorem}
\textit{Proof:}

\hspace{22pt} Let us assume $z$ to be such that $ 0 < \Re(z) < 1$,
$\Im(z) \ne 0$ and

\begin{equation}\label{dois-1}
\zeta(z) = \zeta(1 - z)=0
\end{equation}
and
\begin{equation}\label{dois-2}
\zeta(z) = \zeta(\overline{z})=0
\end{equation}

For \eqref{dois-1} see lemma \ref{lemma-1}. We are using in \eqref{dois-2} the reflection principle
\begin{equation}\label{um-4}
\zeta(\overline{s}) = \overline{\zeta(s)}
\end{equation}

\eqref{dois-1} and \eqref{dois-2} means that the non-trivial zeros lie symmetrically to the real axis
and the line $\Re(s) = \frac{1}{2}$.

\hspace{22pt} The Dirichlet series
\begin{equation}\label{dois-3}
\sum_{n=1}^{\infty} \frac{(-1)^{n-1}}{n^s} =(1 - 2^{1-s})\zeta(s)
\end{equation}

is convergent for all values of $s$ such that $\Re(s) > 0$
~\cite{Titchmarsh}. For $z$ we have
\begin{equation}\label{dois-4}
\zeta(z)= \frac{1}{(1 - 2^{1-z})} \sum_{n=1}^{\infty}
\frac{(-1)^{n-1}}{n^z}
\end{equation}

From \eqref{dois-1}, \eqref{dois-2} and \eqref{dois-4} we obtain
\begin{equation}\label{dois-5}
\frac{1}{(1 - 2^{1 - \overline{z}})} \sum_{n=1}^{\infty}
\frac{(-1)^{n-1}}{n^{\overline{z}}} = \frac{1}{(1 - 2^{z})}
\sum_{n=1}^{\infty} \frac{(-1)^{n-1}}{n^{1-z}}=0
\end{equation}

Since that
\begin{equation}\label{dois-5b}
\frac{1}{(1 - 2^{1 - \overline{z}})} \ne 0
\end{equation}
and
\begin{equation}\label{dois-5c}
 \frac{1}{(1 - 2^{z})} \ne 0
\end{equation}
we have

\begin{equation}\label{dois-6}
\sum_{n=1}^{\infty} \frac{(-1)^{n-1}}{n^{\overline{z}}} =
\sum_{n=1}^{\infty} \frac{(-1)^{n-1}}{n^{1-z}}=0
\end{equation}

Let be
\begin{equation}\label{dois-7}
\Re(z) = \frac{1}{2} + \delta \qquad \textrm{and} \qquad \Im(z) =
t
\end{equation}
where
\begin{equation}\label{dois-8}
- \frac{1}{2} < \delta < \frac{1}{2}
\end{equation}

Substituting \eqref{dois-7} into \eqref{dois-6} we have
\begin{equation}\label{dois-9}
\sum_{n=1}^{\infty} \frac{(-1)^{n-1}}{n^{\frac{1}{2} + \delta -
it}} = \sum_{n=1}^{\infty} \frac{(-1)^{n-1}}{n^{\frac{1}{2} -
\delta - it}}=0
\end{equation}
or
\begin{equation}\label{dois-10}
\sum_{n=1}^{\infty}
\frac{1}{n^{2\delta}}\frac{(-1)^{n-1}}{n^{1-z}} =
\sum_{n=1}^{\infty} \frac{(-1)^{n-1}}{n^{1-z}}=0
\end{equation}

Considering the lemma \ref{lemma-2} we conclude
\begin{equation}\label{dois-11}
\Re(z) = \frac{1}{2}
\end{equation}

\hspace{400pt} $\Box$

\end{document}